\documentclass[runningheads]{llncs}
\usepackage{graphicx}
\usepackage{amsmath, amsfonts, amssymb, amscd, bm, dsfont}

\usepackage[pdftex,
            bookmarks=true,
            bookmarksnumbered=true,
            pdfpagemode=UseOutlines,
            pdfstartview=FitH,
            pdfpagelayout=OneColumn,
            colorlinks=false,
            pdfborder={0 0 0}]{hyperref}
\usepackage{comment}
\newcommand{\real}{\mathbb{R}}

\begin{document}
\title{SHORT-MT: Optimal Solution of Linear Ordinary Differential Equations by Conjugate Gradient Method\thanks{Supported by NSFC 11471307 and CAS Research Program of Frontier Sciences
(QYZDB-SSW-SYS026).}}
\titlerunning{Optimal Solution of Linear Ordinary Differential Equations}

\author{Wenqiang Yang\inst{1},
Wenyuan Wu\inst{1}$\thanks{Corresponding author: \email{wuwenyuan@cigit.ac.cn}.}$, \and Robert M. Corless\inst{2}}
\authorrunning{W. Yang, W. Wu and R. Corless}
\institute{Chongqing Institute of Green and Intelligent Technology, Chinese Academy of Sciences, Chongqing, China \and
Department of Applied Mathematics, Western University, Canada}
\maketitle              % typeset the header of the contribution
\begin{abstract}
Solving initial value problems and boundary value problems of Linear Ordinary Differential Equations (ODEs) plays an important role in many applications. There are various
numerical methods and solvers to obtain approximate solutions represented by points. However, few work about
optimal solution to minimize the residual can be found in
the literatures. In this paper, we first use Hermit cubic spline interpolation at mesh points to represent the solution, then we define the residual error as the $L_2$ norm of the residual obtained by substituting the interpolation solution back to ODEs.
Thus, solving ODEs is reduced to an optimization problem in curtain solution space which can be solved by conjugate gradient method with taking advantages of sparsity of the corresponding matrix.
The examples of IVP and BVP in the paper show that this method can find a solution with smaller global error without additional mesh points.
\end{abstract}
\section{Introduction}\label{sec:intro}

Differential equations (DEs) are one of the most fundamental tools in physical world to model the dynamics of a system.
In machine learning, we may model the learner as some dynamical system, for example a neural network with weights changing according to certain rules. For example, the continuous time recurrent neural network
using a system of ODE to simulate the effects on a neuron of the incoming spike train
works successfully in evolutionary robotics.
Spiking neural networks as third generation of neural network increase the level of realism in a neural simulation
where the neural voltage is usually described as DEs.

Interestingly, people also exploit reversely a variety of neural network methods for solving DEs arising in science and engineering \cite{YNYAKM15}.
In this paper, starting with the most fundamental case: linear ODE,  we attempt to study the optimal solution first by symbolic transformation to an optimization problem and then attain the solution by numerical methods.

The standard from of a Linear ODEs is
\begin{equation}\label{livp}
{\bm x}'={\bm A}(t)\cdot{\bm x}(t)+{\bm q}(t)
\end{equation}
Where ${\bm x}={\bm x}(t):\real\rightarrow\real^n$
is vector of solutions as a function time, the $n\times n$ matrix ${\bm A}(t)$ is the state matrix,
the $n$ dimensional vector ${\bm q}(t)$ corresponds to the inhomogeneous
part of the system.

Such ODE often appearing in numerical modelling and simulation is of great importance in mechanical engineering
and industrial design. Numerical solving of linear ODE is well studied, especially IVP\cite{Walter1998} and BVP\cite{Uri2009}\cite{Shampine2000}\cite{Ascher1988} for ODE. Even a number of efficient solvers
have been developed e.g. Matlab. To qualify and measure the reliability of these solvers without knowing the exact solution is an
important question\cite{Wenyuanwu2018}. A natural choice is to consider the residual error of the equation after substituting the approximate solution\cite{Shampine2005}\cite{Robert2013}\cite{Higham1989}.

Accuracy of the solution depends on the used discrete mesh. If we fix the mesh and consider the solution with first-order continuity approximated by
degree $3$ piecewise polynomials. The key question is to look for a solution with smallest residual error in such solution space which is different from the goal of \cite{Robert2013} to search optimal interpolant with a given numerical solution of ODE.
Our main contribution of this paper is to convert the IVP/BVP of ODEs to a quadratic programming problem. Moreover, the optimal solution can be obtained efficiently by the conjugate gradient method since the corresponding matrix is symmetric and positive definite.

\section{Preliminaries}\label{sec:preli}

\subsection{Hermite Cubic Spline}\label{ssec:hcs}
For first-order continuity, we apply Hermite cubic spline interpolation to construct an approximate solution of ODEs \cite{Erwin2005}.

\begin{theorem}\label{thm:cubicHermit}
 Let $\{t_{0},...,t_{m}\}$ be a set of knots in the interval $[t_{0},t_{m}]$ with
 $t_{0}< t_{1} <...< t_{m}$ and $m\geq2$, there is a unique cubic spline $\widetilde{x} \in C^{1}[t_{0},t_{m}]$ on
 the interval $[t_{i},t_{i+1}]$, for $i = 0,...,m-1$ such that
 \[
    \widetilde{x}(t_{i})= x_{i}, \widetilde{x}'(t_{i}) =x'_{i}
 \]
Let $h_{i}=t_{i+1}-t_{i}$, and $\tau=(t-t_{i})/h_{i}$, $\tau \in [0,1]$, then writing the spline $\widetilde{x}$ as
 \[
    \widetilde{x}_{i}(\tau)=\alpha_{0}(\tau)\cdot x_{i-1}+h_{i}\alpha_{1}(\tau)\cdot x'_{i-1}+\beta_{0}(\tau)\cdot x_{i}+h_{i}\beta_{1}(\tau)\cdot x'_{i}
 \]
 Where $\alpha_{0}(\tau)=2\tau^{3}-3\tau^{2}+1$, $\alpha_{1}(\tau)=\tau^{3}-2\tau^{2}+\tau$, $\beta_{0}(\tau)=-2\tau^{3}+3\tau^{2}$, $\beta_{1}(\tau)=\tau^{3}-\tau^{2}$.
 \end{theorem}

\subsection{Residual and Forward Error}\label{ssec:res}
To evaluate the accuracy of this solution, here we first define the residual and forward error which is introduced in \cite{Wenyuanwu2018}. More detailed study about the error analysis can be found in the book \cite{Robert2013}.

\begin{definition}
Let $n$-dimensional vector $\widetilde{{\bm x}}(t) \in C^1[t_{0},t_{m}]$ be approximate solution of Eq.(\ref{livp}).
Then its \textbf{residual} is defined to be
\begin{equation}\label{re:resdual}
{\bm\delta}(t)=\widetilde{{\bm x}}'(t)-{\bm A}(t)\cdot\widetilde{{\bm x}}(t)-{\bm q}(t)
\end{equation}
The residual of $C^1$ solution is also called \textbf{defect} in \cite{Robert2013}.

The $L_2$ norm of the residual $$\sqrt{\int_{t_0}^{t_{m}} \delta(t)^{T}\delta(t) dt}$$ is called the \textbf{residual error}.

Let $n$-dimensional vector ${\bm x}^{*}(t)$ be the exact solutions of linear ordinary system Eq.(\ref{livp}), such that
\begin{equation}\label{exact solution of lvip}
	({\bm x}^{*})'={\bm A}(t)\cdot{\bm x}^{*}+{\bm q}(t)
\end{equation}
The difference $\Delta{\bm x}(t) = {\bm x}^{*} - \widetilde{{\bm x}}$ is called the \textbf{forward error}.
\end{definition}
In addition, we define the $L_2$ norm of the forward error $$\sqrt{\int_{t_0}^{t_{m}} \Delta{\bm x}^{T}\Delta{\bm x}\; dt}$$ as the \textbf{global error}.

\section{Optimal Solution}\label{sec:opt}
%\subsection{Objective Function}\label{ssec:obf}
In order to minimize the residual error of ODEs solution on the whole time interval, we can define an objective function as follow:
\begin{equation}\label{eq:opt}
\min\sqrt{\sum\limits_{i=0}^{m-1}\int_{t_{i}}^{t_{i+1}} {\bm\delta}_{i}(t)^{T}\cdot {\bm\delta}_{i}(t)\; dt}
\end{equation}
where ${\bm\delta}_{i}(t)$ is the residual on the $i$-th time interval.
Recall Theorem \ref{thm:cubicHermit}, we can express each ${\bm\delta}_{i}(t)$ in terms of ${\bm x}_{i}$, ${\bm x}'_{i}$, ${\bm x}_{i+1}$, ${\bm x}'_{i+1}$ at the mesh points.

%\subsection{Matrixing Method}\label{ssec:mm}
We denote the matrices ${\bm a}_{i}=\frac{\alpha'_{0}}{h_{i}}{\bm E}-\alpha_{0}{\bm A}$, ${\bm b}_{i}=\alpha'_{1}{\bm E}-h_{i}\alpha_{1}{\bm A}$, ${\bm c}_{i}=\frac{\beta'_{0}}{h_{i}}{\bm E}-\beta_{0}{\bm A}$ and ${\bm d}_{i}=\beta'_{1}{\bm E}-h_{i}\beta_{1}{\bm A}$, where ${\bm E}$ is the $n\times n$ identity matrix. Thus, $${\bm\delta}_{i}={\bm a}_{i}{\bm x}_{i}+{\bm b}_{i}{\bm x}'_{i}+{\bm c}_{i}{\bm x}_{i+1}+{\bm d}_{i}{\bm x}'_{i+1}-{\bm q}.$$

Additionally, let $I_{i}=\int_{t_{i}}^{t_{i+1}} {\bm\delta}_{i}(t)^{T}\cdot{\bm\delta}_{i}(t)\; dt$.
Then we have
\begin{equation}
I_{i}={\bm y}^{T}_{i}{\bm F}_{i}{\bm y}_{i}-2{\bm B}_{i}^{T}{\bm y}_{i}+\int_{t_{i}}^{t_{i+1}}{\bm q}^{T}{\bm q}dt
\end{equation}
where  ${\bm F}_{i}=\int_{t_{i}}^{t_{i+1}}\left( \begin {array}{cccc}
                    {\bm a}^{T}_{i}{\bm a}_{i}&{\bm a}^{T}_{i}{\bm b}_{i}&{\bm a}^{T}_{i}{\bm c}_{i}&{\bm a}^{T}_{i}{\bm d}_{i}\\
                      {\bm b}^{T}_{i}{\bm a}_{i}& {\bm b}^{T}_{i}{\bm b}_{i}&{\bm b}^{T}_{i}{\bm c}_{i}&{\bm b}^{T}_{i}{\bm d}_{i}\\
                     {\bm c}^{T}_{i}{\bm a}_{i}&{\bm c}^{T}_{i}{\bm b}_{i}&{\bm c}^{T}_{i}{\bm c}_{i}&{\bm c}^{T}_{i}{\bm d}_{i}\\
                     {\bm d}^{T}_{i}{\bm a}_{i}&{\bm d}^{T}_{i}{\bm b}_{i}&{\bm d}^{T}_{i}{\bm c}_{i}&{\bm d}^{T}_{i}{\bm d}_{i}
                   \end {array} \right)dt\;\;$ is a $4n\times 4n$ matrix, and  \\
                    ${\bm B}_{i}=\int_{t_{i}}^{t_{i+1}}{\bm q}^{T}\left( \begin {array}{c}
                    {\bm a}_{i}\\
                    {\bm b}_{i}\\
                    {\bm c}_{i}\\
                    {\bm d}_{i}
                    \end {array} \right)dt\;\;$ and   ${\bm y}_{i}=\left( \begin {array}{c}
                    {\bm x}_{i}\\
                    {\bm x}'_{i}\\
                    {\bm x}_{i+1}\\
                    {\bm x}'_{i+1}
                    \end {array} \right)$ are  $4n$-dimensional vectors.

Here ${\bm x}_{i}$, ${\bm x}'_{i}$, ${\bm x}_{i+1}$, ${\bm x}'_{i+1}$ are unknowns in the optimization problem Eq.(\ref{eq:opt}).

\begin{figure}
  \includegraphics[height=6cm, width=10cm]{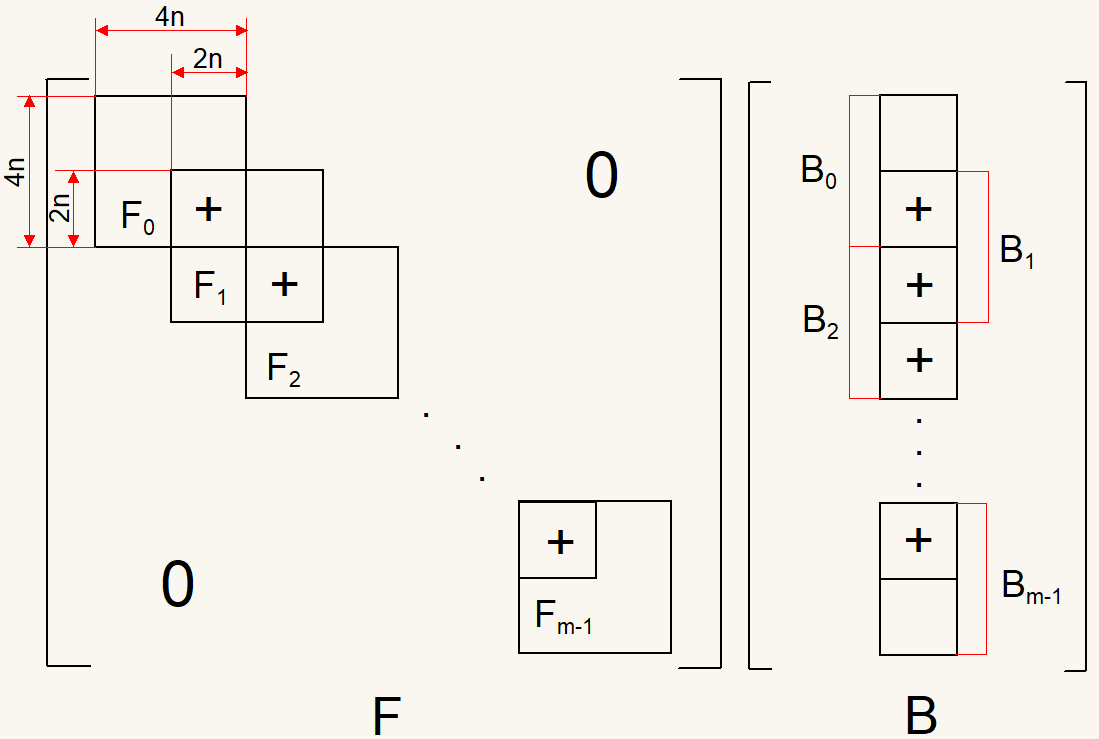}\\
  \caption{The Structure of Matrix F and Vector B. Here ``$+$" indicates the overlapping of two adjacent blocks.}
  \label{fig:structure}
\end{figure}

%\subsection{Conjugate Gradient Method}\label{ssec:cgm}
Since $\sum\limits_{i=0}^{m-1}\int_{t_{i}}^{t_{i+1}}{\bm q}^{T}{\bm q}dt\geq 0$ is a constant, to minimize $I_{i}$ is equivalent to minimize ${\bm y}^{T}_{i}{\bm F}_{i}{\bm y}_{i}-2{\bm q}^{T}_{i}{\bm B}_{i}{\bm y}_{i}$. It's obviously that objective function can be rewritten as a quadratic form
\begin{equation}\label{eq:quad}
opt=\min{\bm y}^{T}{\bm F}{\bm y}-2{\bm B}^{T}{\bm y}_{i}.
\end{equation}
where ${\bm F}$ is highly structured with size $2n(m+1) \times 2n(m+1)$ (see Fig. \ref{fig:structure}).

It is not difficult to verify that matrix ${\bm F}$ is a symmetric positive definite matrix. According to \cite{ref_JR Shewchuk}, the minimal value can be attained at the solution of
\begin{equation}
{\bm F}{\bm y}={\bm B}
\end{equation}
 It is well-known that such linear system can be solved efficiently by using conjugate gradient method. Especially, if we have obtained an approximate numerical solution by other solvers, we can use it as an initial point of conjugate gradient method to refine the solution.

It is important to point out here that the same formulation Eq.(\ref{eq:quad}) works for both initial value problems (IVPs) and boundary value problems (BVPs) of Linear ODEs. The only difference is that the first $n$ rows and columns of ${\bm F}$ will be removed after substituting the initial values for IVPs. But for BVPs, there are $n$ rows in different places and the corresponding columns which will be removed. Accordingly, the substitution also happens in ${\bm B}$.  Thus, the number of unknowns of Eq.(\ref{eq:quad}) drops by $n$.

\section{Examples}\label{sec:exam}
In order to show the performance of our method to IVPs and BVPs for linear ODEs, we give an RC ladder network system as an example, which can be used to filter a signal by blocking certain frequencies and passing others. When we consider its time response, this system can be considered as an IVP. In addition, when we consider its controllability, this system becomes a BVP.
We use a constant state matrix $\textbf{A}$ since it is easy to obtain the exact solution for comparison.

\begin{figure}
  \center{
  \includegraphics[height=3cm, width=7cm]{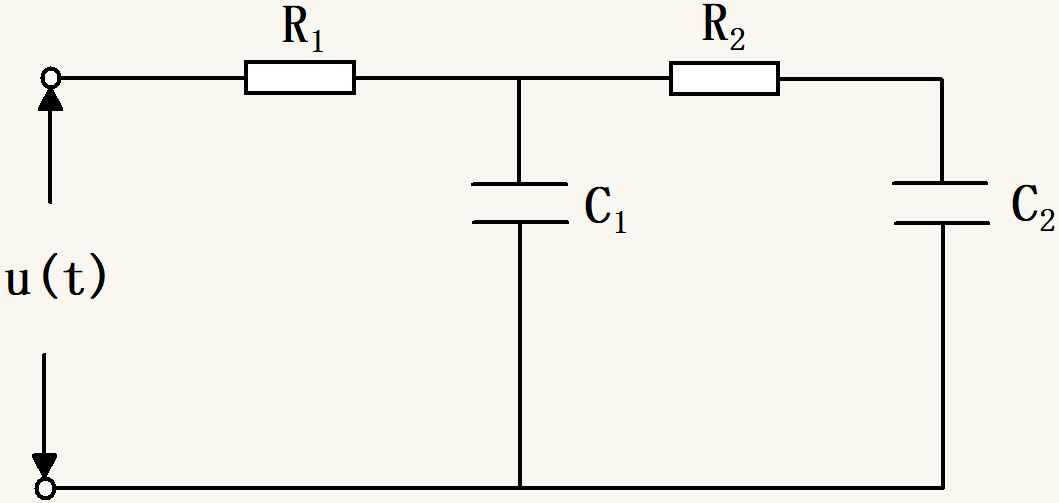}\\
  \caption{An RC Ladder Network}
  \label{fig:rccircle}
  }
\end{figure}

\subsection{An Initial Value Problem}\label{ssec:ivp}
An RC ladder network is shown in Fig. \ref{fig:rccircle}. Assume that the state variables are ${\bm x}=[x_{1},x_{2}]^{T}$ be the voltage across each capacitor. During the time period $t\in[0,2]$, we study the time response function of the system. Suppose the resistances $R_{1}=R_{2}=10k\Omega$, and the capacitance $C_{1}=C_{2}=100\mu F$. At the beginning, ${\bm x}(0)=[0,0]^{T}$, we give an input signal $u(t)=2\cdot sin(2t)$.
So the model can be described as follows:
\[
{\bm x}'=\left( \begin {array}{rr}
                    -2& 1\\1&-1
                   \end {array} \right){\bm x}+\left( \begin {array}{c}
                    1\\0
                   \end {array} \right)u(t)
\]

\begin{figure}
  \includegraphics[height=4cm, width=12cm]{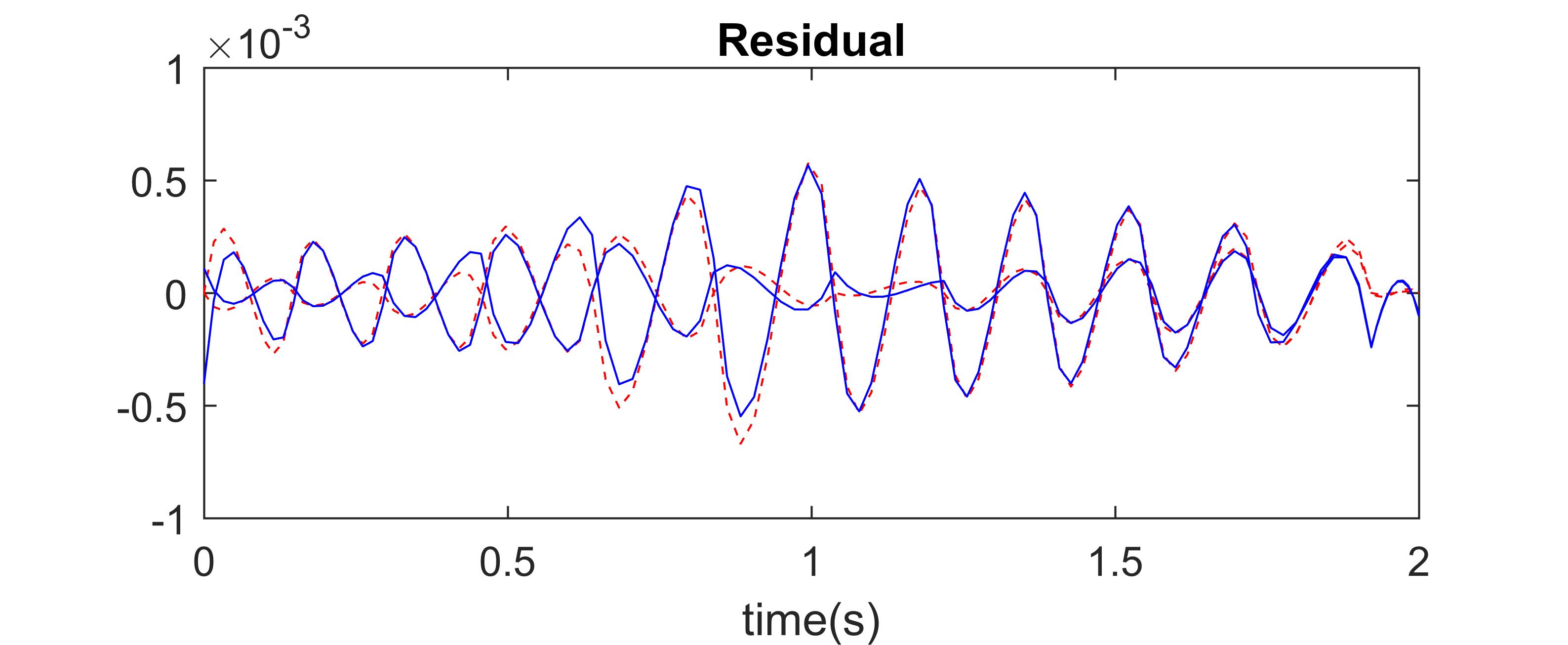}\\
  \caption{Residual comparison for IVP between ODE45 (red dash) and our method (blue). }
  \label{fig:ivp_res}
\end{figure}

\begin{figure}
  \includegraphics[height=8cm, width=12cm]{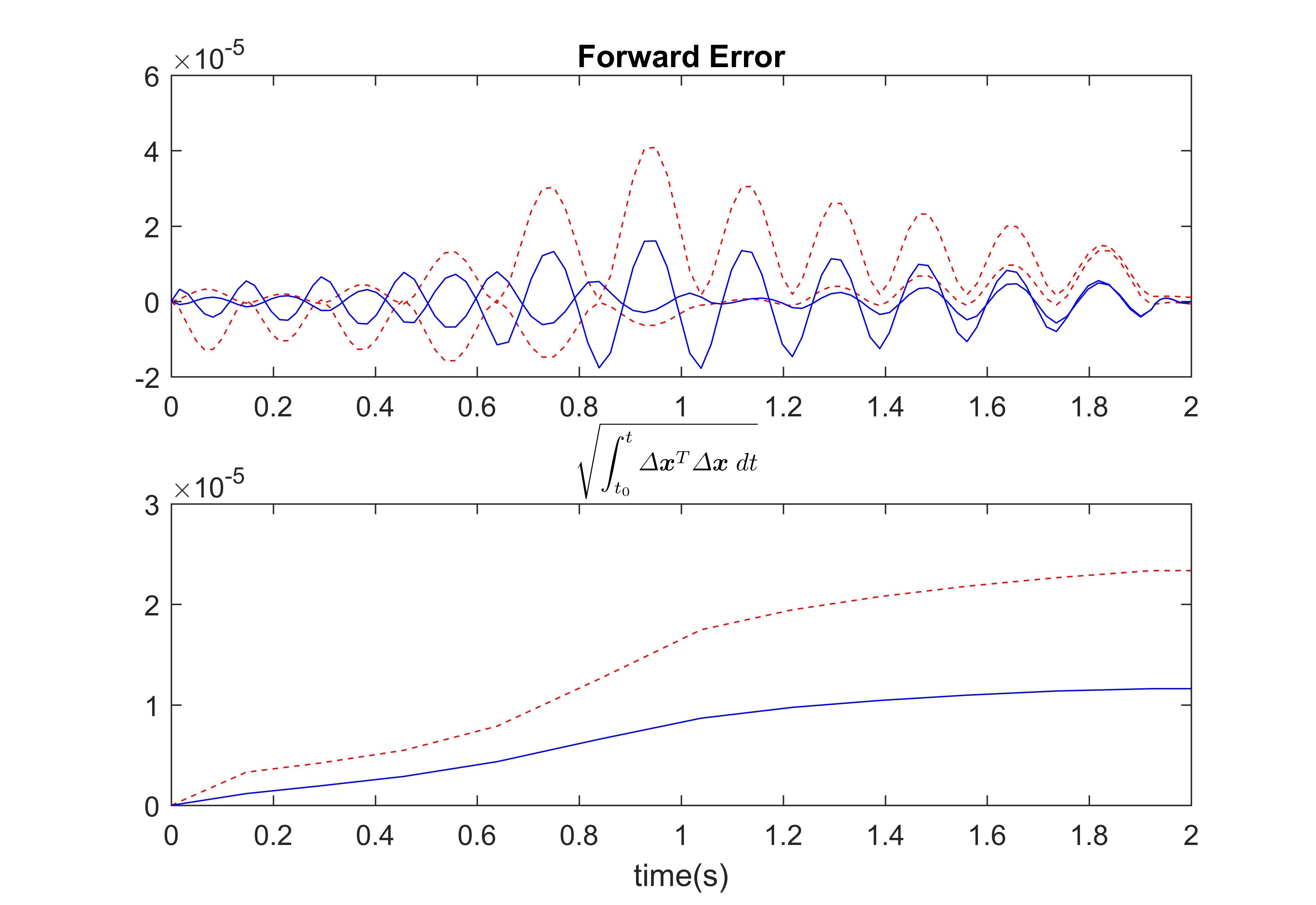}\\
  \caption{Error comparison for IVP between ODE45 (red dash) and our method (blue). }
  \label{fig:ivp}
\end{figure}

Obviously, this ODE has an exact solution as ${\bm x}^{*}=e^{{\bm A}t}\cdot \int_{0}^{t}e^{-{\bm A}s}{\bm q}(s)ds$. We apply \textit{ode$45$} and our method, which is implemented also in Matlab, to the same $13$ mesh-points decided by \textit{ode$45$}, and get the approximate solutions respectively. The corresponding $F$ is a $50\times50$ matrix and its condition number is $3.949\times 10^{3}$. The residual of both methods are almost the same in Fig.\ref{fig:ivp_res}, and residual error of our method  $0.4097\times 10^{-3}$ is slightly less than $0.4275\times 10^{-3}$ by \textit{ode$45$}. However, Fig. \ref{fig:ivp} shows that global error by our method (blue curve) $9.1313\times 10^{-5}$ is much smaller than the one by \textit{ode$45$} (red dash curve)$1.8232\times 10^{-4}$.  It's clear that the forward errors of both methods are fluctuating around zero but the blue curve given by our method has a smaller amplitude than \textit{ode$45$}. The reason behind would be our method is a kind of global approach compared with the local approach \textit{ode$45$}.

\subsection{A Boundary Value Problem}\label{ssec:bvp}
If we want to control the former RC ladder network to achieve the goal $x_{1}(0)=0$ and $x_{2}(2)=1$, and input signal $u(t)=2\cdot sin(2t)$, how about the beginning voltage of $x_{2}(0)$? To answer the question, we need to solve bvp for ODEs.
\begin{figure}
  \includegraphics[height=4cm, width=12cm]{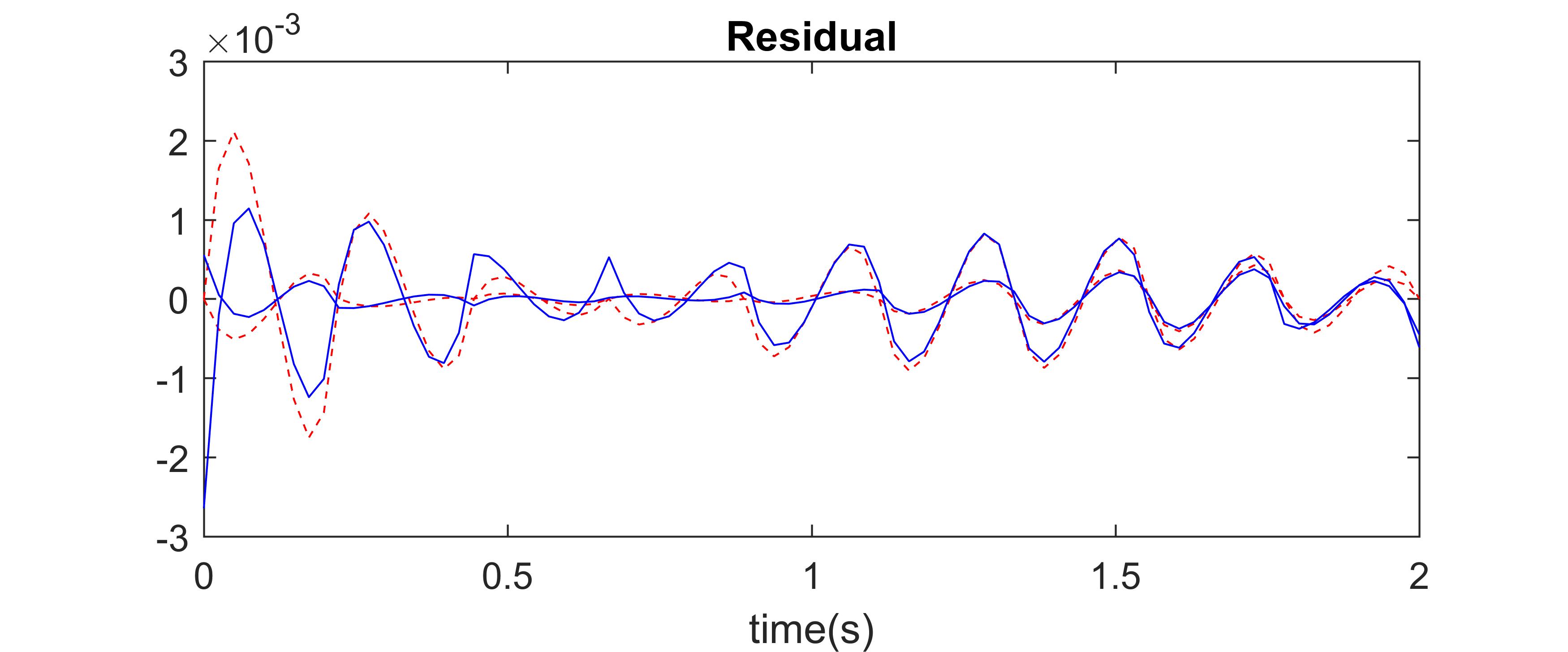}\\
  \caption{Residual comparison for IVP between ODE45 (red dash) and our method (blue). }
  \label{fig:bvp_res}
\end{figure}

\begin{figure}
  \includegraphics[height=8cm, width=12cm]{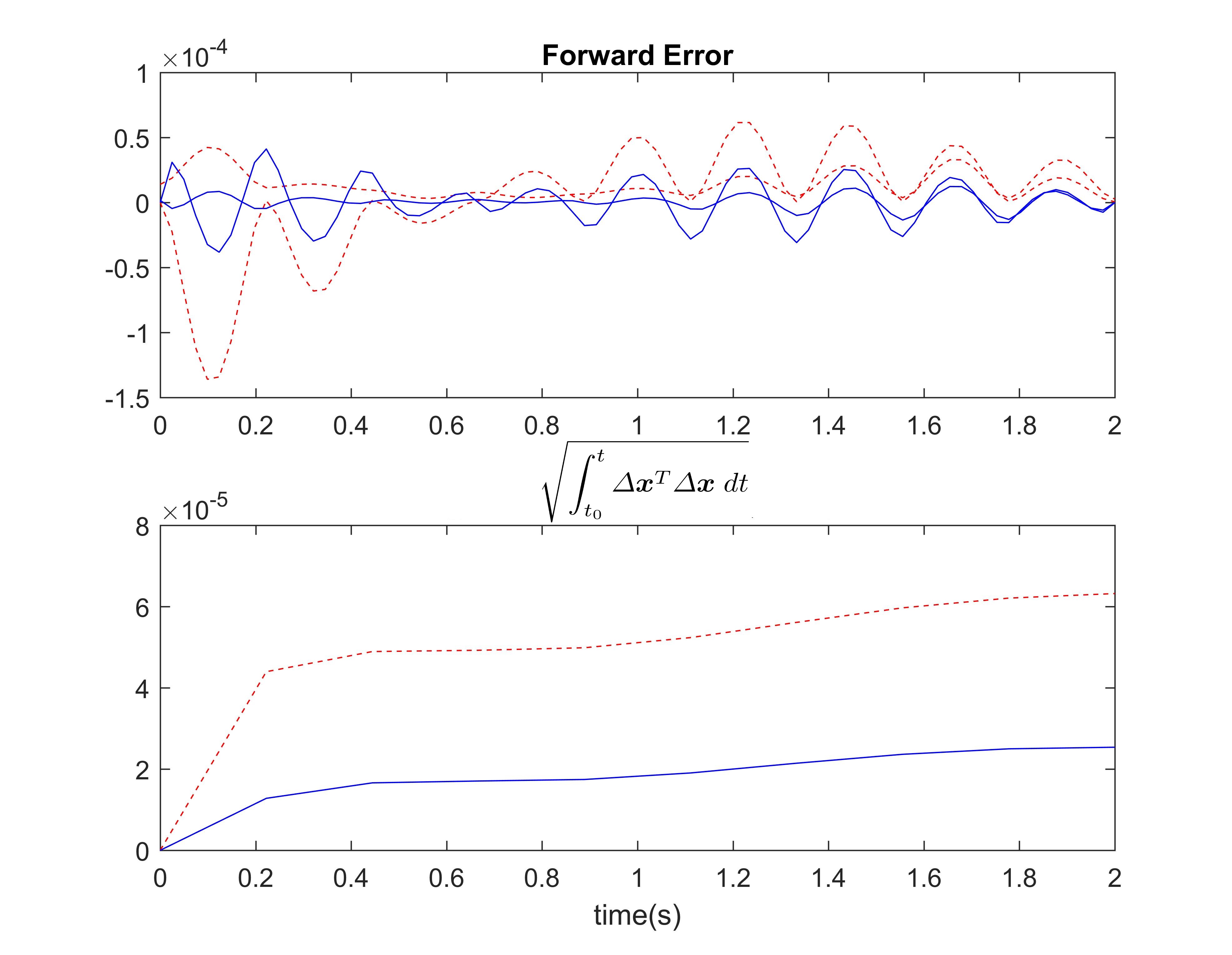}\\
  \caption{Error comparison for BVP between ODE45 (red dash) and our method (blue).}
  \label{fig:bvp}
\end{figure}

In this example, We apply Matlab's \textit{\textit{bvp5c}} and our method to the same $10$ mesh-points decided by \textit{\textit{bvp5c}}. In our method, $F$ is a $38\times38$ matrix with the condition number $1.640\times 10^{3}$. In order to compare the forward error, we also compute its exact solution by Sec 6.1 \cite{Uri2009}. The residual of both methods are almost the same in Fig.\ref{fig:bvp_res}, and residual error of our method $0.7942\times 10^{-3}$ is less than $0.9595\times 10^{-3}$ given by \textit{\textit{bvp5c}}. Similar to the IVP case, the same comparison result can be observed in Fig. \ref{fig:bvp}, that global error by our method (blue curve) $1.7824\times 10^{-4}$ is much smaller than the one by \textit{\textit{bvp5c}} (red dash curve) $4.8521\times 10^{-4}$.

\section{Conclusion}\label{sec:con}
In this paper, we provide a conjugate gradient based method for computing optimal solution of linear ordinary differential equations. A fine mesh in time interval leads to a large but sparse matrix $\textbf{F}$. Due to nice properties of $\textbf{F}$
the iteration approach for linear system can be applied here and it can take advantages of sparsity.
This method can solve both IVPs and BVPs of Linear ODEs.

As an illustrative example, we use RC ladder network to generate IVP and BVP of Linear ODEs, whose exact solution can be found easily. By comparison, we find that our method can give solutions with smaller forward error and global error than solvers provided by Matlab at the same mesh-points.

\end{document}